\documentclass[11pt]{article}

\usepackage{amsthm,amsfonts,amssymb,amsmath,mathrsfs}
\usepackage{diagrams}
\usepackage{graphicx,color}
\usepackage{srcltx}
\usepackage[latin1]{inputenc}  
\usepackage[T1]{fontenc}
\usepackage[all,dvips]{xy}
\newarrow{Dots}{}...{}
\newarrow{Tosurj}{}---{->>}
\newarrow{Toinj}{boldhook}---{->}
\newarrow{Map}{}---{->}
\newarrow{Isom}{boldhook}---{->>}
\newarrow{Double}{}==={}

\def\hw{{\hat w}}
\def\hG{{\hat G}}
\def\hg{{\hat g}}

\def\hT{{\hat T}}
\def\hB{{\hat B}}
\def\hlambda{{\hat\lambda}}
\def\hW{{\hat W}}
\def\halpha{{\hat \alpha}}\def\hbeta{{\hat \beta}}

\newcounter{paragrafsubsub}[subsubsection]
\renewcommand{\theparagrafsubsub}{%
\thesubsubsection.\roman{paragrafsubsub}}
\newcommand{\paragrafsubsub}{%
\refstepcounter{paragrafsubsub}
{\bf \theparagrafsubsub}\hspace{0.2em}--- }

\newcounter{paragrafsub}[subsection]
\renewcommand{\theparagrafsub}{\thesubsection.\arabic{paragrafsub}}
\newcommand{\paragrafsub}{%
\refstepcounter{paragrafsub}
{\bf \theparagrafsub}\hspace{0.2em}--- }

\newcounter{paragraf}[section]
\renewcommand{\theparagraf}{\thesection.\arabic{paragraf}}
\newcommand{\paragraf}{%
\refstepcounter{paragraf}
{\bf \theparagraf}\hspace{0.2em}--- }

\newcommand\paragraphe{%
\par \indent
\ifcase\value{subsection} %
\paragraf
\else
\ifcase\value{subsubsection}\paragrafsub %
\else\paragrafsubsub
\fi\fi
}

\def\longto{\longrightarrow}

\def\lr{{\mathcal{LR}}}

\def\PP{{\mathbb P}}
\def\QQ{{\mathbb Q}}\def\ZZ{{\mathbb Z}}
\def\CC{{\mathbb C}}

\def\Gl{{\rm Gl}}
\def\Spin{{\rm Spin}}
\def\SO{{\rm SO}}

\def\Pic{\rm Pic}

\def\SL{{\rm SL}}\def\Sp{{\rm Sp}}
\def\SO{{\rm SO}}

\def\Li{{\mathcal{L}}}
\def\Mi{{\mathcal{M}}}
\def\O{{\mathcal O}}

\def\LR{{\rm LR}}

\def\GL{{\rm Gl}}
\def\SP{{\rm Sp}}
\def\h{\hat}

\newtheorem{lemma}{Lemma}
\newtheorem{prop}{Proposition}
\newtheorem{theo}{Theorem}

\newenvironment{remark}{{\noindent\bf Remark.}}{}

\begin{document}
\title{Two generalizations of the PRV conjecture}
\author{PL. Montagard\footnote{Universit{\'e} Montpellier II - 
CC 51-Place Eug{\`e}ne Bataillon -
34095 Montpellier Cedex 5 -
France - {\tt pierre-louis.montagard@math.univ-montp2.fr}},
B. Pasquier\footnote{Universit{\'e} Montpellier II - 
CC 51-Place Eug{\`e}ne Bataillon -
34095 Montpellier Cedex 5 -
France - {\tt boris.pasquier@math.univ-montp2.fr}},
N. Ressayre\footnote{Universit{\'e} Montpellier II - 
CC 51-Place Eug{\`e}ne Bataillon -
34095 Montpellier Cedex 5 -
France - {\tt nicolas.ressayre@math.univ-montp2.fr}. The author was partially supported by the French National Research Agency (ANR-09-JCJC-0102-01).}
}

\maketitle

\begin{abstract}
Let $G$ be a complex connected reductive group.
The PRV conjecture, which  was proved 
independently by S.~Kumar and O.~Mathieu in 1989, gives explicit irreducible 
submodules of the tensor product of two irreducible $G$-modules.
This paper has three aims. First, we simplify the proof of the PRV conjecture,
then we generalize it to other branching problems.
Finally, we find other irreducible components of  the tensor product of two irreducible $G$-modules
that appear for ``the same reason'' as the PRV ones.

\end{abstract}

\section{Introduction}

\subsection{The original PRV conjecture}

Parthasarathy-Ranga Rao-Varadarajan conjectured in the sixties the following

\bigskip
\noindent{\it The PRV conjecture.}
Let $G$ be a complex connected reductive group  with associated Weyl group $W$.
Let $V_G(\mu)$ and $V_G(\nu)$ be two irreducible  $G$-modules with highest weights
$\mu$ and $\nu$ respectively.
Then, for any $w\in W$, the irreducible $G$-module $V_G(\overline{\mu+w\nu})$ 
with extremal weight
$\mu+w\nu$, occurs with multiplicity at least one in $V_G(\mu)\otimes V_G(\nu)$.\\

This conjecture was proved independently by S.~Kumar in \cite{Kumar:prv1} and
O.~Mathieu in \cite{Mathieu:prv}.
The aim of this paper is to simplify the proof of the PRV conjecture and to
generalize it in two directions.

\subsection{Two generalizations}

We now assume that $G$ is a subgroup of a bigger connected reductive group $\hat G$.
Fix a Borel subgroup $\hat B$ and a maximal torus $\hat T\subset \hat B$ 
of $\hat G$ such that
$B=\hat B\cap G$ is a Borel subgroup of $G$ and $T=\hat T\cap G$ is a maximal torus
 of $G$.
Consider the restriction map $\rho\,:\,X(\hat T)\longto X(T)$ from the character
group of $ \hat T$ to the one of $T$.
Let $\hat \lambda$ be a dominant weight of $\hat T$ and $V_{\hat G}(\hat \lambda)$ be the irreducible 
$\hat G$-module of highest weight $\hat \lambda$. Let $\hat w\in\hat W$. 
The first aim of this paper is the following

\bigskip\noindent
{\bf Question.} Does the irreducible $G$-module 
$V_G(\overline{\rho(\hat w\hat \lambda)})$ with
extremal weight $\rho(\hat w\hat\lambda)$ occur  with multiplicity at least one in 
$V_{\hat G}(\hat \lambda)$?

\bigskip
Although the answer may be NO (examples are given in
Section~\ref{sec:SL3SO3} or in Section~\ref{sec:Kro}), the PRV 
conjecture exactly asserts that the answer is YES if $G$ is diagonally embedded in 
$\hat G=G\times G$.

Let $\hat G/\hat B$ denote the complete flag variety of $\hat G$, 
$X_\hw^\circ$ denote the $G$-orbit
$G\hw\hB/\hB$ and $X_\hw$ denote its closure in $\hG/\hB$. If $X_\hw^\circ$ is closed in $\hG/\hB$, we easily
check that the answer is YES. 
We also answer positively  the question under a topological assumption on $X_\hw$.

\begin{theo}
\label{thi:ssgrpe}
We assume $X_\hw$ is multiplicity free.

Then, $V_G(\overline{\rho(\hw\hat\lambda)})$ is a $G$-submodule of $V_\hG(\hat\lambda)$.
\end{theo}

Here, $X_\hw$ is said to be multiplicity free if its cycle class in the cohomology of
$\hG/\hB$ is a linear combination with coefficients 0 or 1 of Schubert classes.
This assumption, which can be hard to check, is fulfilled for example 
if $G$ is a spherical subgroup
of $\hG$ of minimal rank (see \cite{spherangmin} for the complete list of such subgroup).
In particular, $G$ is a spherical subgroup of $G\times G$ of minimal rank and Theorem~\ref{thi:ssgrpe} implies the original PRV conjecture.

\bigskip
Our second generalization of the PRV conjecture deals with the decomposition
 of  tensor products: we exhibit new components. 

\begin{theo}
\label{thi:tensorprod}
Let $\lambda,\,\mu,\,\nu$ be three dominant weights of $T$.
We assume that there exist $v,\,w\in W$, a simple root $\alpha$ and  an integer $k$ 
such that 
\begin{enumerate}
\item \label{ass:long}$l(s_\alpha v)=l(v)+1$, $l(s_\alpha w)=l(w)+1$; 
\item \label{ass:rellmn}$\lambda=v\mu+w\nu-k\alpha$;
\item \label{ass:ineq}$0\leq k\leq\langle v\mu,\alpha^\vee\rangle$, and
$0\leq k\leq\langle w\nu,\alpha^\vee\rangle$.
\end{enumerate}
Then, $V_G(\lambda)$ is a submodule of $V_G(\mu)\otimes V_G(\nu)$.
\end{theo}

Here, $\alpha^\vee$ denotes the coroot associated to $\alpha$, and
$\langle\cdot,\cdot\rangle$ denotes the pairing between weights and coroots.
We obtain the original PRV conjecture by applying Theorem~\ref{thi:tensorprod} with extremal values
of $k$ in \ref{ass:ineq}. 

\subsection{About proofs}

The two key ingredients in our proofs are the normality of $X_\hw$, and the fact that for any
$\hG$-linearized and globally generated line bundle $\Li$ on $\hG/\hB$, the restriction map 
${\rm H}^0(\hG/\hB,\Li)\longto  {\rm H}^0(X_\hw,\Li)$
 is surjective (see~Theorem~\ref{th:ssgrpe2} below). 
 An analogue version of these two results was already proved by M.~Demazure in the case of any Schubert varieties in flag varieties \cite{Dem:norm&surj}, before S.~Kumar proved them in the case where $\hG=G\times G$ \cite{Kumar:prv1}. In our context, we need to use the generalization of these results for any $G$, $\hG$ and multiplicity free $X_{\hat w}$, obtained by M.~Brion in \cite{Br:multfree}.
These two ingredients also play a central role in Kumar's proof. 
But, Kumar's proof also uses a complete description of ${\rm H}^0(X_\hw,\Li)$ 
mainly due to Bott and the Joseph filtration. We make these two latter ingredients useless
by using an argument of semistability.

\subsection{Link with a saturation problem}

In the general situation 
 $G\subset\hat G$, we consider 
the set $\LR(G,\hG)$ of pairs $(\lambda,\hlambda)$ of dominant weights of $T$ 
and $\hat T$ such that $V_G(\lambda)$ occurs in $V_\hG(\hlambda)$.
By a Brion-Knop's theorem, $\LR(G,\hG)$ is a finitely generated semigroup.
From a theoretic viewpoint, the convex cone $\lr(G,\hG)$ 
generated by $\LR(G,\hG)$ is well understood: the complete and minimal list 
of inequalities  is parametrized 
by explicit cohomological conditions (see~\cite{GITEigen}).
There are so many inequalities that it is not  obvious to concretely 
describe this cone and especially to construct points in this cone.
A starting point in the proof of Theorem~\ref{thi:ssgrpe} is the following 

\begin{prop}\label{prop:asymp}
  Let $\hlambda$ be a dominant character  of $\hT$ and $\hw\in \hW$.

Then, there exists a positive integer $n$ such that 
$V_G(n\overline{\rho(\hw\hlambda)})$ is a $G$-submodule of $V_\hG(n\hlambda)$.

In other words, $(\overline{\rho(\hat w\hat\lambda)},\hat\lambda)$ belongs to 
$\lr(G,\hG)$.
\end{prop}

With the additional assumption that $X_\hw$ is multiplicity free, Theorem~\ref{thi:ssgrpe}
asserts that  $(\overline{\rho(\hat w\hat\lambda)},\hat\lambda)$ belongs to 
$\LR(G,\hG)$.
The question of understanding the difference between $\lr(G,\hG)$ and $\LR(G,\hG)$ is
known as a saturation problem.
Let $\Lambda$ be the subgroup of $X(T)\times X(\hT)$ 
generated by $\LR(G,\hG)$. The semigroup $\LR(G,\hG)$ is said to be {\it saturated along a half line} if the first non-zero point of $\Lambda$ on this half line belongs to $\LR(G,\hG)$ (and $\LR(G,\hG)$ is said to be {\it saturated} if it is along any half line in $\lr(G,\hG)$).  
 Theorem~\ref{thi:ssgrpe} shows that if $X_\hw$ is multiplicity free, $\LR(G,\hG)$ is saturated along all the half lines  given by Proposition~\ref{prop:asymp}.

Knutson-Tao proved in \cite{KT:saturation} that $\LR(\SL_n,\SL_n\times\SL_n)$ is saturated.
Belkale-Kumar proved in \cite{BK:typeBC} that $\LR(\Sp_{2n},\Sp_{2n}\times\Sp_{2n})$ 
and  $\LR(\Spin_{2n+1},\Spin_{2n+1}\times\Spin_{2n+1})$ are saturated up to a factor 2:
the second point of $\Lambda$ on any half-line belongs to $\LR$.
Kapovich-Leeb-Millson obtained important results on the saturation question for 
semigroups $\LR(G,G\times G)$ (see \cite{KLM:tri}).

We can now explain Theorem~\ref{thi:tensorprod} in this context. 
Fix two dominant weights $\mu$ and $\nu$ of $T$. 
The intersection of $\lr(G,G\times G)$ with $X(T)\otimes\QQ\times\{\mu\}\times\{\nu\}$
is a polytope $P(\mu,\nu)$ (namely, a moment polytope). 
The original PRV conjecture gives finitely many  points in $P(\mu,\nu)$ that 
generate saturated half lines.
Theorem~\ref{thi:tensorprod}  gives finitely many segments in $P(\mu,\nu)$
whose all rational points generate saturated half lines
(see Section~\ref{sec:expletensor} for examples).

\section{Restriction to a subgroup}
\label{sec:ssgrpe}

\subsection{Setting}

Let $G$ be a complex connected reductive group, 
with a fixed Borel subgroup $B$ and maximal torus $T\subset B$.
Let $X(T)$ denote the character group of $T$.
For any dominant weight $\lambda\in X(T)$, let $V_G(\lambda)$ denote the irreducible
$G$-module with highest weight $\lambda$.
Let $W$ be the Weyl group of $(G,T)$. For any character $\lambda$, 
the orbit $W.\lambda$ intersects the dominant chamber in one point
denoted by $\bar\lambda$.
We will denote by $w_0$ the longest element of the Weyl group $W$.\\

We now assume that $G$ is a subgroup of a connected   reductive group  $\hG$. 
Let $\hT$ and $\hB$ be a maximal torus and a Borel subgroup of $\hG$ such that 
$T\subset\hT\subset\hB\supset B$.
We will use hats to denote objects relative to $\hat G$ instead of $G$; for example we
will write $\hat W$, $\hat w_0$, \dots\ 
For a given dominant character $\hlambda$ of $\hT$, we are interested in the following\\

\noindent{\sc Problem.} Find irreducible $G$-submodules of $V_\hG(\hlambda)$?

\subsection{$G$-orbits in the complete flag manifold of $\hG$}

For any $\hw\in \hW$, we set $X_\hw^\circ=G\hw\hB/\hB$ and $X_\hw$ its closure.
We also denote by $\sigma_\hw$ the cycle class of the Schubert variety $\overline{\hB \hw\hB/\hB}$
in $\hG/\hB$.
It is well known that
\begin{eqnarray}
  \label{eq:cohomGB}
  {\rm H}^*(\hG/\hB,\ZZ)=\bigoplus_{\hw\in \hW}\ZZ.\sigma_\hw.
\end{eqnarray}

Let $V$ be an irreducible subvariety of $\hG/\hB$.
The cycle class $[V]$ of $V$ in ${\rm H}^*(\hG/\hB,\ZZ)$ can be expanded as follows
\begin{eqnarray}
  \label{eq:Vsigma}
  [V]=\sum_{\hw\in \hW}c_\hw(V)\sigma_\hw,
\end{eqnarray}
where the $c_\hw(V)$ are non-negative integers.
The variety $V$ is said to be multiplicity free if for any $\hw\in \hW$, $c_\hw(V)=0$ or $1$.

\subsection{The statement}
Consider the restriction map
$
\rho\,:\,X(\hT)\longto X(T).
$
We now state a slightly more general version of Theorem~\ref{thi:ssgrpe}.

\begin{theo}
\label{th:ssgrpe}
With above notation, let $\hlambda$ be a dominant character  of $\hT$ and $\hw\in \hW$.
We assume that one of the following assumption holds:
\begin{enumerate}
\item\label{hyp:fermee}
 $X_\hw^\circ$ is closed;
\item \label{hyp:spherangmin}
$G$ is spherical of minimal rank in $\hG$;
\item \label{hyp:multfree}
$X_\hw$ is multiplicity free;
\item \label{hyp:multfreedual}
$X_{\hw \hw_0}$ is multiplicity free.
\end{enumerate}

Then, $V_G(\overline{\rho(\hw\hlambda)})$ is a $G$-submodule of $V_\hG(\hlambda)$.
\end{theo}

The first case is easy and certainly well known.

\begin{proof}[Proof in case~\ref{hyp:fermee}]
Since $X_\hw^\circ$ is complete, the isotropy group of $\hw\hB/\hB$ in $G$ is a parabolic 
subgroup of $G$.
But, it is contained in $\hw\hB\hw^{-1}$, so it is solvable. 
It follows that $B':=\hw\hB\hw^{-1}\cap G$ is a Borel subgroup of $G$ containing $T$.
Then there exists $w\in W$ such that $w^{-1}Bw=B'$.

Let $v$ be a non-zero vector of $V_\hG(\hlambda)$ of highest weight $\hlambda$. 
It is clear that $\hw v$ is an eigenvector of weight $\rho(\hw\hlambda)$ for $B'$ (here, we
identify $X(T)$ and $X(B')$ by the restriction morphism). 
It follows that $w\hw v$ is an eigenvector of weight $w\rho(\hw\hlambda)$ for $B$, 
so that $w\rho(\hw\hlambda)$ is dominant and $w\rho(\hw\hlambda)=\overline{\rho(\hw\hlambda)}$.
The theorem follows.
\end{proof}

We now prove case~\ref{hyp:multfreedual} assuming that 
case~\ref{hyp:multfree} is known.

\begin{proof}[Proof in  case~\ref{hyp:multfreedual}]
We apply the theorem in case~\ref{hyp:multfree} to the 
dominant weight $-\hw_0\hlambda$ and the element $\hw \hw_0$ of $\hW$. 
We obtain that $V_G(\overline{\rho(-\hat w \lambda)})$ is contained in
$V_{\hat G}(-\hat w_0\hlambda)=V_\hG(\hlambda)^*$.
Since $-\rho(\hat w \lambda)$ is an extremal weight of
$V_G(\overline{\rho(\hat w \lambda)})^*$,
we deduce that   $V_G(\overline{\rho(\hw\hlambda)})^*$ is a $G$-submodule of $V_\hG(\hlambda)^*$.
The theorem follows by duality.
\end{proof}

\subsection{The  spherical case}

In this subsection, we look at a situation where we can check when
assumption~\ref{hyp:multfree} is fulfilled. This will allow us to
discuss the various assumptions on examples and to include
case~\ref{hyp:spherangmin}
in case~\ref{hyp:multfree}.

\bigskip
\paragraphe
Assume that $G$ is a spherical subgroup of $\hG$; 
i.e. $G$ acts on $\hG/\hB$ with finitely many orbits.
In \cite{Br:GammaGH}, M.~Brion defined an oriented graph $\Gamma(\hG/G)$ 
whose vertices are the $G$-orbit closures in $\hG/\hB$. The edges, which can be simple or double, are labeled by the simple roots of $\hG$.
The assumption ``$X_\hw$ is multiplicity free'' can be easily read off this graph:
$X_\hw$ is multiplicity free if and only if for any path from $X_\hw$ to $\hG/\hB$ there is no double
edge.
In particular, by \cite[Proposition~2.1]{spherangmin}, 
if $G$ is spherical of minimal rank, any $G$-orbit 
closure in $\hG/\hB$ is multiplicity free.
In particular, case~\ref{hyp:spherangmin} of Theorem~\ref{th:ssgrpe} is a consequence of
case~\ref{hyp:multfree}.\\

We now study two examples where $G$ is spherical, which illustrate Theorem~\ref{th:ssgrpe}.

\bigskip
\paragraphe
Let $\hat G=\Sp_4$ and $G=\Gl_2$ be the Levi subgroup of a maximal parabolic subgroup of $\Sp_4$ that stabilizes an isotropic plane in $\CC^4$. Then $G$ is a spherical subgroup of $\hat G$ and the oriented graph $\Gamma(\hat G/G)$ (with arrows pointed down) is the following ($\halpha$ and $\hbeta$ denote respectively the short and the long  simple roots of $\Sp_4$).
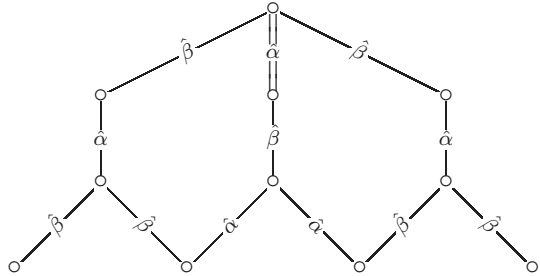
\begin{diagram}[abut,height=1.5em,width=1.5em]
&&&&&&\circ&&&&&\\
&&&&&\ldLine(4,2)~\hbeta&\dDouble~\halpha&\rdLine(4,2)~\hbeta&&&&\\
&&\circ&&&&\circ&&&&\circ&&\\
&&\dLine~\halpha&&&&\dLine~\hbeta&&&&\dLine~\halpha&\\
&&\circ&&&&\circ&&&&\circ&&\\
&\ldLine~\hbeta&&\rdLine~\hbeta&&\ldLine~\halpha&&\rdLine~\halpha&&
\ldLine~\hbeta&&\rdLine~\hbeta\\
\circ&&&&\circ&&&&\circ&&&&\circ\\
\end{diagram}
In this example, the varieties $X_{\hat w}$ correspond to the four vertices at the bottom of the graph $\Gamma(\hat G/G)$ and they are in fact the four closed $G$-orbits in $\hG/G$. So Theorem~\ref{th:ssgrpe} can be applied here for all $\hat w\in \hat W$. This gives an example where we need to use hypothesis~\ref{hyp:fermee} of Theorem~\ref{th:ssgrpe} to apply it, because two of the closed $G$-orbits above are not multiplicity free.\\

\paragraphe\label{sec:SL3SO3}
Let $\hat G=\SL_3$ and $G=\SO_3$ naturally embedded in $\SL_3$. 
Let $\alpha$, $\halpha$ and $\hbeta$ denote the simple roots of $\SO_3$ and $\SL_3$.
Also denote by $\omega_\alpha$, $\omega_\halpha$ and $\omega_\hbeta$ the corresponding fundamental 
weights.
Then $G$ is a spherical subgroup of $\hat G$ and the oriented graph $\Gamma(\hat G/G)$ 
is the following diagram.
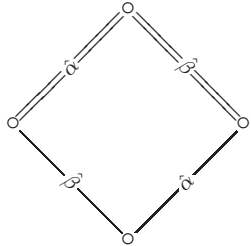
\begin{diagram}[abut,height=2em,width=2em]
&&\circ&&&&\\
&\ldDouble~\halpha&&\rdDouble~\hbeta\\
\circ&&&&\circ&\\
&\rdLine~\hbeta&&\ldLine~\halpha\\
&&\circ&&\\
\end{diagram}
  
We can read on the graph that there exist exactly two not-closed $G$-orbits  in $\hG/\hB$ with multiplicity, 
namely $X_{s_\halpha}^\circ$ and $X_{s_\hbeta}^\circ$. 
An easy computation gives us that $\rho(s_\halpha\omega_\halpha)=\rho(s_\hbeta\omega_\hbeta)=0$. 
But we can also check that $V_G(0)$ is neither in $V_{\hat G}(\omega_\halpha)$ nor in $V_{\hat G}(\omega_\hbeta)$, so that Theorem~\ref{th:ssgrpe} is not satisfied in these two cases. 

We have just seen that $(0,\omega_\halpha)$ is not in the semi-group $\LR(G,\hG)$ defined in the introduction. But we can remark that $(0,2\omega_\halpha)\in \LR(G,\hG)$, while $(0,\omega_\halpha)$ is in the subgroup of $X(T)\times X(\hT)$ generated by $\LR(G,\hG)$ (because we can compute that $(2\omega_\alpha,\omega_\hbeta)$ and $(2\omega_\alpha,\omega_\halpha+\omega_\hbeta)$ are in $\LR(G,\hG)$). Then $\LR(G,\hG)$ is not saturated along the half line generated by $(0,2\omega_\halpha)$.

\bigskip
The rest of Section~\ref{sec:ssgrpe} is devoted to the proof of
Theorem~\ref{th:ssgrpe}.

\subsection{A result of Geometric Invariant Theory}

Let $X$ be  an irreducible projective  $G$-variety.
As in \cite{GIT}, we denote by  $\Pic^G(X)$ the group of $G$-linearized line bundles on $X$.
Let  $\Li\in\Pic^G(X)$ and  let ${\rm H}^0(X,\Li)$ denote the $G$-module of regular sections 
of $\Li$.
A point $x\in X$ is said to be {\it semistable with respect to $\Li$}  if there exists $ n>0$
 and 
$\tau\in{\rm H}^0(X,\Li^{\otimes n})^G$ such that $\tau(x)\neq 0$.

\bigskip
\begin{remark}
  Note that this definition of semistable points is not standard.
Indeed, it is usually agreed that the open
subset defined by the non-vanishing of $\tau$ is affine. This property, which is useful to 
construct a good quotient, is automatic only if $\Li$ is ample; 
hence, 
our definition coincides with the usual one if $\Li$ is ample.
\end{remark}\\

A line bundle $\Li$ on $X$ is said to be {\it semiample} if a positive power of $\Li$ 
is base point free.
 If $\Li$ is a line bundle on $X$ and $x$ is a point in $X$, $\Li_x$ denotes 
the fiber in $\Li$ over $x$. 
We will need the following lemma mainly due to Kostant.

\begin{lemma}
\label{lem:Kos}
Let $\Li\in\Pic^G(X)$ be semiample and  $x\in X$ be a $T$-fixed point. 
We assume that $T$ acts
trivially on $\Li_x$.

Then $x$ is semistable with respect to $\Li$.
\end{lemma}

\begin{proof}
  Let $n$ be a positive integer, such that the natural morphism
$$
\varphi\,:\,X\longto\PP({\rm H}^0(X,\Li^{\otimes n})^*)
$$
is well defined. Set $V={\rm H}^0(X,\Li^{\otimes n})^*$.
Let $v\in V$ be a non-zero vector on the line $\varphi(x)$.
The assumption implies that $v$ is fixed by $T$.

Let $U$ be the unipotent radical of $B$. Then, as an orbit of an unipotent group in
an affine variety, $U.v$ is closed in $V$ (see \cite[Theorem~2]{Ros:unip}); and, 
$B.v=U.v$. Since $G/B$ is complete, it follows that $G.v$ is closed in $V$.
We deduce that there exists a $G$-invariant homogeneous polynomial $P$ of degree $d$ 
on $V$ such that $P(\varphi(x))\neq 0$. It follows that there exists a $G$-invariant
section $\tau$ of $\Li^{\otimes nd}$ such that $\tau(x)\neq 0$.
\end{proof}

\subsection{The Borel-Weil theorem}

Let $P$ be a parabolic subgroup of $G$.
Let $\nu$ be a character of $P$. Let $\CC_\nu$ denote the field $\CC$ 
endowed with the action of $P$ defined by
$p.\tau=\nu(p^{-1})\tau$ for all $\tau\in\CC_\nu$ and $p\in P$. We define the fiber product 
$G\times_P \CC_\nu$ as the quotient of $G\times\CC_\nu$ by the following equivalent relation $$\forall g\in G,\,\forall \tau\in\CC_\nu\mbox{ and }\forall p\in P,\,\,(g,\tau)\sim(gp,p^{-1}.\tau).$$ 
It is a $G$-linearized 
line bundle on $G/P$, denoted  by $\Li_\nu$. 
In fact, the map $$\begin{array}{ccc}X(P)&\longto&\Pic^G(G/P)\\ \nu&\longmapsto&\Li_\nu\end{array}$$ is an isomorphism.

We assume that $P$ contains  $B$ (in that case, $P$ is said to be {\it standard}).
Then,  $X(P)$ identifies with a subgroup of $X(T)$.
For $\nu\in X(P)$, $\Li_\nu$ is semiample if and only if it has non-zero sections
if and only if $\nu$ is dominant.
Moreover, ${\rm H}^0(G/P,\Li_\nu)$ is the dual of the simple $G$-module of highest weight $\nu$.
For $\nu$ dominant, $\Li_\nu$ is ample if and only if $\nu$ cannot be extended to a subgroup of $G$
bigger than $P$.\\

\subsection{The Brion theorem}

We will need the following theorem, due to Brion, on multiplicity free subvarieties of $G/B$.

\begin{theo}\cite[Theorem~1]{Br:multfree}
  \label{th:Brion}
Let $V$ be a multiplicity free subvariety of $G/B$ and $\Li$ be any semiample $G$-linearized line bundle on $G/B$,  then
\begin{enumerate}
\item $V$ is normal;
\item the restriction map
${\rm H}^0(G/B,\Li)\longto {\rm H}^0(V,\Li)$
is surjective.
\end{enumerate}
\end{theo}

\subsection{Proof of Theorem~\ref{th:ssgrpe}}

\paragraphe
We first prove an {\bf asymptotic version} of  Theorem~\ref{th:ssgrpe}, that is 
Proposition~\ref{prop:asymp} of the introduction.



\begin{proof}[Proof of Proposition~\ref{prop:asymp}]
Set $X=\hG/\hB$. By the Borel-Weil theorem, we have
$$
{\rm H}^0(X,\Li_{\hlambda})=V_\hG(\hlambda)^*.
$$
It remains to prove that, for some $n>0$, $\Li_{\hlambda}^{\otimes n}$ admits a non-zero section that is an 
eigenvector of weight $-n\overline{\rho(\hw\hlambda)}$ for the opposite Borel subgroup $B^-$ of $G$.
This is made more precisely in Lemma~\ref{lem:asym} below.
\end{proof}

\begin{lemma}
 \label{lem:asym} 
There exists $n>0$ such that 
$\Li_{\hlambda}^{\otimes n}$ admits a section $\tau$ which is an 
eigenvector of weight $-n\overline{\rho(\hw\hlambda)}$ for $B^-$ such that
the restriction of $\tau$ to $X_\hw^\circ$ is non-zero.
\end{lemma}

\begin{proof}
Consider the variety $Y=X\times G/B^-$ endowed with the diagonal action of $G$ given by:
$g'.(\hg\hB/\hB,gB^-/B^-)=(g'\hg\hB/\hB,g'gB^-/B^-)$.
Let $\Li_{-\overline{\rho(\hw\hlambda)}}^-$ be the $G$-linearized line bundle on $G/B^-$ such that
$B^-$ acts on the fiber over $B^-$ by the character $\overline{\rho(\hw\hlambda)}$.
We also consider the line bundle $\Mi:=\Li_\hlambda\boxtimes\Li^-_{-\overline{\rho(\hw\hlambda)}}$ on $Y$. Note that $\Mi$ is semi-ample because $\hlambda$ is dominant and $-\overline{\rho(\hw\hlambda)}$ is dominant with respect to $B^-$.

By definition of $\overline{\rho(\hw\hlambda)}$, there exists $v\in W$ such that 
$\overline{\rho(\hw\hlambda)}=v\rho(\hw\hlambda)$. 
Then, it is clear that $T$ acts trivially on the fiber in $\Mi$ over the 
point $y:=(v\hw\hB/\hB,B^-/B^-)$. 
Now, applying Lemma~\ref{lem:Kos}, we obtain, for some $n>0$, a section $\tau_Y\in H^0(Y,\Mi^{\otimes n})^G$ such that $\tau_Y(y)\neq 0$.

Define $\tau$ as the restriction of $\tau_Y$ to $X\times B^-/B^-$ seen as a section  
of $\Li_\hlambda$ on $X$. 
Since $\tau_Y$ is $G$-invariant, $\tau$ is $B^-$-equivariant of weight $-n\overline{\rho(\hw\hlambda)}$. 
 It is clear that $\tau(v\hw\hB/\hB)\neq 0$, so that the restriction of $\tau$  to $X_\hw^\circ$ is non-zero. The lemma is proved.
\end{proof}

\paragraphe
We have already seen that it is sufficient to
prove Theorem~\ref{th:ssgrpe} under assumption~\ref{hyp:multfree}.
By Theorem~\ref{th:Brion}, it is sufficient to prove the following

\begin{theo}
\label{th:ssgrpe2}
Let $\hlambda$ be a dominant character  of $\hT$ and $\hw\in \hW$.
We assume that 
\begin{enumerate}
\item\label{hyp:normal}
 $X_\hw$ is normal; 
\item \label{hyp:surj}
 the restriction map
$
{\rm H}^0(\hG/\hB,\Li_{\hlambda})\longto {\rm H}^0(X_\hw,\Li_{\hlambda})
$
is surjective.
\end{enumerate}

Then, $V_G(\overline{\rho(\hw\hlambda)})$ is a $G$-submodule of $V_\hG(\hlambda)$.
\end{theo}

\begin{proof}
Consider the following restriction maps:
$$\xymatrix{
{\rm H}^0(\hG/\hB,\Li_{\hlambda})\ar@{->>}[r]&
{\rm H}^0(X_\hw,\Li_{\hlambda})\ar[r]&
{\rm H}^0(X_\hw^\circ,\Li_{\hlambda}).\\}
$$
Since the first one is surjective and $G$-equivariant, 
it is sufficient to find $V_G(\overline{\rho(\hw\hlambda)})^*$ in 
${\rm H}^0(X_\hw,\Li_{\hlambda})$.
We will first prove that 
$V_G(\overline{\rho(\hw\hlambda)})^*$ is a submodule of 
${\rm H}^0(X_\hw^\circ,\Li_{\hlambda})$ without multiplicity. 
Next, we will prove that the corresponding 
$B^-$-equivariant section on $X_\hw^\circ$ extends to $X_\hw$ using both the asymptotic version and 
the normality of $X_\hw$.
  
\bigskip
By Lemma~\ref{lem:H0X0} below, there exists a (unique up to scalar multiplication) 
non-zero regular section $\sigma$ of $\Li_\hlambda$ on
$X_\hw^\circ$ which is $B^-$-equivariant of weight $\overline{-\rho(\hw\hlambda)}$.

Let $n>0$ and $\tau$ be as in Lemma~\ref{lem:asym}.
Then, $\tau_{|X_\hw^\circ}$ and $\sigma^{\otimes n}$ are two non-zero regular sections
 of $\Li_{n\hlambda}$ on
$X_\hw^\circ$ which are $B^-$-equivariant of weight $-n\overline{\rho(\hw\hlambda)}$.
By Lemma~\ref{lem:H0X0}, it follows that $\tau_{|X_\hw^\circ}$ and $\sigma^{\otimes n}$ are proportional.
In particular, $\sigma^{\otimes n}$ extends to a regular section of $\Li_{n\hlambda}$ on $X_\hw$.
Since $X_\hw$ is normal, it follows that $\sigma$ also extends to a regular section of 
$\Li_{\hlambda}$ on $X_\hw$.
The theorem is proved.
\end{proof}

\noindent{\bf Notation.}
If  $H$ is an algebraic affine group, $\chi$ is a character of $H$ and $V$ is a $H$-module, we set:
$$
V^{(H)_\chi}=\{v\in V\mid\,\forall h\in H,\ h.v=\chi(h)v\}.
$$

\begin{lemma}
\label{lem:H0X0}
  The $G$-module $V_G(\overline{\rho(\hw\hlambda)})^*$ has multiplicity exactly one in 
${\rm H}^0(X^\circ_\hw,\Li_{\hlambda})$.
\end{lemma}

\begin{proof}
Let $G_\hw\subset G$ be the isotropy group of $\hw\hB/\hB$ so that $X_\hw^\circ$ is
isomorphic to the homogeneous space $G/G_\hw$. 
Let us define $\mu=\rho({\hw\hlambda})$.
Since $G_\hw$ acts on the
fiber $(\Li_{\hlambda})_{\hw\hB/\hB}$ by the character $-\mu$, the
line bundle $\Li$ on $X_\hw^\circ$ is isomorphic to $G\times_{G_\hw}\CC_{-\mu}$.

Then by using the Frobenius decomposition, the space of global sections 
$H^0(G/G_\hw,G\times_{G_\hw}\CC_{-\mu})$ can be identify with: 

$$\bigoplus_{\chi}V^*_G(\chi)\otimes
  (V_G(\chi))^{(G_\hw)_{\mu}},$$ 
where the sum is over the set of dominant weights of $G$.  
So, we have to prove that the vector space $V_G(\overline{\mu})^{(G_\hw)_{\mu}}$ is one-dimensional. First, since $G_\hw=G\cap \hw\hB\hw^{-1}$ contains $T$,
 the dimension of $V_G(\overline{\mu})^{(G_\hw)_{\mu}}$ is smaller than one.
  
  The dimension is exactly one if $G_\hw$ is contained in  the
  parabolic group $P_G(\mu)$ associated to the weight $\mu$.
 By Lemma~\ref{lem:asym}, there exist an integer $n$ and a section $\tau\in
  H^0(X_\hw,\Li_{n\hlambda})^{(B^-)_{n\overline{\mu}}}$ such that the
  restriction of $\tau$ to $X^\circ_\hw$
  is non-zero. So the dimension of
  $H^0(X^\circ_\hw,\Li_{n\hlambda})^{(B^-)_{n\overline{\mu}}}$ is bigger than
  one. By using the Frobenius decomposition as above, we deduce that the dimension of
  $V_G(\overline{n\mu})^{(G_\hw)_{n\mu}}$ is bigger than (and so equal
  to) one, and that the parabolic group $P_G(n\mu)$ associated to
  the weight $n\mu$ contains the group $G_\hw$. We conclude by saying that $P_G(n\mu)=P_G(\mu)$. 
\end{proof}

\section{Applications}

\subsection{Applications to the Kronecker product}\label{sec:Kro}

The aim of this section is to detail our results for
$\GL(E)\times\GL(F)\subset\GL(E\otimes F)$.
This problem is equivalent to the  question on the
decomposition of tensor products of representations for the symmetric group.

A partition $\pi$ is a sequence
$\pi=(\pi_1,\pi_2,\ldots,\pi_k)$ of weakly decreasing non-negative
integers. 
By convention, we allow partitions with some zero parts, and two partitions
that differ by zero parts are the same. If several parts are equal we denote the multiplicity of this part by an exponent. For
example $(3^2,2^4,1)$ means the partition $(3,3,2,2,2,2,1)$.
For any  partition $\pi$, we define 
$|\pi|=\pi_1+\pi_2+\cdots+\pi_k$ and $l(\pi)$ as the number of non-zero parts of $\pi$.

Recall that if $V$ is
a finite dimensional vector space, then the 
 $\GL(V)$-irreducible polynomial representations  are
 in bijection with the partitions $\pi$ such 
that $l(\pi)\leq \dim V$: we denote by $S_\pi V$ the representation
associated to $\pi$.

Let $E,F$ be two vector spaces of respective dimension $m,n$, and
consider $G=\GL(E)\times \GL(F)$ and $\h G=\GL(E\otimes F)$. 
Let $\gamma$ be a partition such that $l(\gamma)\leq mn$. We can decompose the
irreducible
representation $S_\gamma(E\otimes F)$ as a $G$-representation: 
$$S_\gamma(E\otimes F)=\sum_{\alpha,\beta}
N_{\alpha\beta\gamma}\,S_\alpha E\otimes S_\beta F\,,$$
where the sum is taken over partitions $\alpha,\beta$ such that
$|\alpha|=|\beta |=|\gamma|$, $l(\alpha)\leq m$ and $l(\beta)\leq n$. 

\bigskip
\begin{remark}
The irreducible representations of the symmetric group  $\mathfrak S_n$ 
correspond bijectively with the partitions $\pi$ such that $|\pi|=n$; we denote by
$[\pi]$ the  representation corresponding to $\pi$.
By  using the Schur-Weyl duality, we can show that $N_{\alpha\beta\gamma}$ is also the
multiplicity of $[\gamma]$ in $[\alpha]\otimes[\beta]$ (see for example
\cite[Chapter~6]{book:FH}). 
Now, the fact that the representations of $\mathfrak S_n$ are
self-dual implies that  $N_{\alpha\beta\gamma}$ is symmetric in
$\alpha$, $\beta$ and $\gamma$.
\end{remark}

\bigskip
By fixing  basis of $E$ and $F$, we denote by $T_E$ and $T_F$ the maximal
tori of $\GL(E)$ and $\GL(F)$ consisting of diagonal matrices.
For $i=1,\ldots, m$,  denote by
$\eta_i$ the character that maps an element of $T_E$
 to its $i^{\rm th}$ diagonal coefficient.
Similarly, we define the characters $\delta_j$'s of $T_F$.
The basis of $E$ and $F$ induce a natural basis of $E\otimes F$ indexed
by pairs $(i,j)$. Let $\hat T$ denote the corresponding
maximal torus of $\hat G$ and  $\h\varepsilon_{i,j}$ the character of
$\h T$ corresponding to $(i,j)$.
Note that $\rho(\h\varepsilon_{i,j})=(\eta_i,\delta_j)$. 

The coordinates of characters of $\h T$ in the basis
$\h\varepsilon_{i,j}$, which are indexed by pairs $(i,j)$, will be 
represented in tableaux of $m$ lines and $n$ columns.  
For any tableau $t$ (identified with the corresponding character of $\h
T$),  $\rho(t)$ is obtained by summing along columns, to obtain the coordinates of a
character of $T_E$, and along lines, to obtain the coordinates of a character of $T_F$.

In Theorem~\ref{th:ssgrpe}, the weights of the form $\h w\h\lambda$
are exactly the extremal weights of $V_G(\h\lambda)$. 
In particular, they do not depend on the choice of a Borel subgroup of
$\h G$ but only on $\h T$ and the representation $V_{\h
  G}(\h\lambda)$.
 Here, we have fixed the torus  and the representation:
the extremal weights of $\h T$ in $S_\gamma(E\otimes F)$ 
are the tableaux $m\times n$ filled by the parts of $\gamma$.

For example, suppose that $m=n=3$ and the two following tableaux
 correspond to extremal weights of $S_{1^4}(E\otimes F)$:

\begin{center}
\setlength{\unitlength}{3947sp}%
\begingroup\makeatletter\ifx\SetFigFont\undefined%
\gdef\SetFigFont#1#2#3#4#5{%
  \reset@font\fontsize{#1}{#2pt}%
  \fontfamily{#3}\fontseries{#4}\fontshape{#5}%
  \selectfont}%
\fi\endgroup%
\begin{picture}(3027,1276)(3364,-2525)
\thinlines
{\color[rgb]{0,0,0}\put(5176,-1261){\line( 1, 0){900}}
\put(6076,-1261){\line( 0,-1){900}}
\put(6076,-2161){\line(-1, 0){900}}
\put(5176,-2161){\line( 0, 1){900}}
}%
{\color[rgb]{0,0,0}\put(5176,-1561){\line( 1, 0){900}}
}%
{\color[rgb]{0,0,0}\put(5176,-1861){\line( 1, 0){900}}
}%
{\color[rgb]{0,0,0}\put(5476,-1261){\line( 0,-1){900}}
}%
{\color[rgb]{0,0,0}\put(5776,-1261){\line( 0,-1){900}}
}%
\put(5626,-2461){\makebox(0,0)[b]{\smash{{\SetFigFont{12}{14.4}{\rmdefault}{\mddefault}{\updefault}{\color[rgb]{0,0,0}$1$}%
}}}}
\put(5926,-2461){\makebox(0,0)[b]{\smash{{\SetFigFont{12}{14.4}{\rmdefault}{\mddefault}{\updefault}{\color[rgb]{0,0,0}$2$}%
}}}}
\put(5326,-2461){\makebox(0,0)[b]{\smash{{\SetFigFont{12}{14.4}{\rmdefault}{\mddefault}{\updefault}{\color[rgb]{0,0,0}$1$}%
}}}}
\put(5326,-1486){\makebox(0,0)[b]{\smash{{\SetFigFont{12}{14.4}{\rmdefault}{\mddefault}{\updefault}{\color[rgb]{0,0,0}$1$}%
}}}}
\put(5926,-1486){\makebox(0,0)[b]{\smash{{\SetFigFont{12}{14.4}{\rmdefault}{\mddefault}{\updefault}{\color[rgb]{0,0,0}$1$}%
}}}}
\put(5926,-2086){\makebox(0,0)[b]{\smash{{\SetFigFont{12}{14.4}{\rmdefault}{\mddefault}{\updefault}{\color[rgb]{0,0,0}$1$}%
}}}}
\put(6376,-1486){\makebox(0,0)[b]{\smash{{\SetFigFont{12}{14.4}{\rmdefault}{\mddefault}{\updefault}{\color[rgb]{0,0,0}$2$}%
}}}}
\put(6376,-1786){\makebox(0,0)[b]{\smash{{\SetFigFont{12}{14.4}{\rmdefault}{\mddefault}{\updefault}{\color[rgb]{0,0,0}$1$}%
}}}}
\put(6376,-2086){\makebox(0,0)[b]{\smash{{\SetFigFont{12}{14.4}{\rmdefault}{\mddefault}{\updefault}{\color[rgb]{0,0,0}$1$}%
}}}}
\put(5626,-1786){\makebox(0,0)[b]{\smash{{\SetFigFont{12}{14.4}{\rmdefault}{\mddefault}{\updefault}{\color[rgb]{0,0,0}$1$}%
}}}}
{\color[rgb]{0,0,0}\put(3376,-1261){\line( 1, 0){900}}
\put(4276,-1261){\line( 0,-1){900}}
\put(4276,-2161){\line(-1, 0){900}}
\put(3376,-2161){\line( 0, 1){900}}
}%
{\color[rgb]{0,0,0}\put(3376,-1561){\line( 1, 0){900}}
}%
{\color[rgb]{0,0,0}\put(3376,-1861){\line( 1, 0){900}}
}%
{\color[rgb]{0,0,0}\put(3676,-1261){\line( 0,-1){900}}
}%
{\color[rgb]{0,0,0}\put(3976,-1261){\line( 0,-1){900}}
}%
\put(3826,-2461){\makebox(0,0)[b]{\smash{{\SetFigFont{12}{14.4}{\rmdefault}{\mddefault}{\updefault}{\color[rgb]{0,0,0}$1$}%
}}}}
\put(3526,-1486){\makebox(0,0)[b]{\smash{{\SetFigFont{12}{14.4}{\rmdefault}{\mddefault}{\updefault}{\color[rgb]{0,0,0}$1$}%
}}}}
\put(4126,-1486){\makebox(0,0)[b]{\smash{{\SetFigFont{12}{14.4}{\rmdefault}{\mddefault}{\updefault}{\color[rgb]{0,0,0}$1$}%
}}}}
\put(3826,-1486){\makebox(0,0)[b]{\smash{{\SetFigFont{12}{14.4}{\rmdefault}{\mddefault}{\updefault}{\color[rgb]{0,0,0}$1$}%
}}}}
\put(4576,-1486){\makebox(0,0)[b]{\smash{{\SetFigFont{12}{14.4}{\rmdefault}{\mddefault}{\updefault}{\color[rgb]{0,0,0}$3$}%
}}}}
\put(3526,-1786){\makebox(0,0)[b]{\smash{{\SetFigFont{12}{14.4}{\rmdefault}{\mddefault}{\updefault}{\color[rgb]{0,0,0}$1$}%
}}}}
\put(3526,-2461){\makebox(0,0)[b]{\smash{{\SetFigFont{12}{14.4}{\rmdefault}{\mddefault}{\updefault}{\color[rgb]{0,0,0}$2$}%
}}}}
\put(4126,-2461){\makebox(0,0)[b]{\smash{{\SetFigFont{12}{14.4}{\rmdefault}{\mddefault}{\updefault}{\color[rgb]{0,0,0}$1$}%
}}}}
\put(4576,-1786){\makebox(0,0)[b]{\smash{{\SetFigFont{12}{14.4}{\rmdefault}{\mddefault}{\updefault}{\color[rgb]{0,0,0}$1$}%
}}}}
\put(4576,-2086){\makebox(0,0)[b]{\smash{{\SetFigFont{12}{14.4}{\rmdefault}{\mddefault}{\updefault}{\color[rgb]{0,0,0}$0$}%
}}}}
\end{picture}%

\end{center}
where the boxes corresponding to zero coordinates are left empty.  

In the first tableau, $\rho(t)=\overline{\rho(t)}=((3,1,0),(2,1,1))$.
We can easily check that  the irreducible
representation  $[1^4]$ (which is the one dimensional
representation given by the signature of $\mathfrak S_4$) appears in
the tensor product $[3,1]\otimes[2,1^2]$.

In the second tableau, $\rho(t)=((2,1,1),(1,1,2))$
and 
$\overline{\rho(t)}=((2,1,1),(2,1,1))$.
We can check that   $[2^4]$ appears in $[4,2^2]\otimes[4,2^2]$ which matches with our
asymptotic result (Proposition~\ref{prop:asymp}).
But, the irreducible
representation  $[1^4]$ does not appear in
the tensor product $[2,1^2]\otimes[2,1^2]$. 
Then, Theorem~\ref{th:ssgrpe} shows that some
$\GL(E)\times\GL(F)$-orbit closures of the form $X_{\h w}$ of the  complete flag variety of $E\times F$
are not multiplicity free. 
A natural but probably difficult question appears here:
which orbit closures  $X_{\h w}$ (for $\h w\in \h W$) are multiplicity free?

\bigskip
We can prove that the $\h w$'s in $\h W$ such that  the orbit $X_{\h w}^\circ$ is
closed, correspond bijectively to standard tableaux $m\times n$.
Now case~\ref{hyp:fermee} of Theorem~\ref{th:ssgrpe} 
gives the following rule to compute some components of the
tensor product of two representations of the symmetric group. We don't
know if this rule is already known.

\bigskip
\noindent
{\bf Rule.} {\it 1. Fill the tableau $m\times n$ by the parts of $\gamma$
in weakly decreasing order along lines and columns.\\
2. Sum along lines and columns to obtain $\alpha$ and $\beta$.

Then, $[\gamma]$ appears in $[\alpha]\otimes[\beta]$.}
\bigskip

For example, the tableaux: 

\vspace*{.5cm}
\begin{center}
\setlength{\unitlength}{3947sp}%
\begingroup\makeatletter\ifx\SetFigFont\undefined%
\gdef\SetFigFont#1#2#3#4#5{%
  \reset@font\fontsize{#1}{#2pt}%
  \fontfamily{#3}\fontseries{#4}\fontshape{#5}%
  \selectfont}%
\fi\endgroup%
\begin{picture}(5427,1351)(-1136,-2600)
\thinlines
{\color[rgb]{0,0,0}\put(-1124,-1561){\line( 1, 0){900}}
}%
{\color[rgb]{0,0,0}\put(-1124,-1861){\line( 1, 0){900}}
}%
{\color[rgb]{0,0,0}\put(-1124,-2161){\line( 1, 0){900}}
}%
{\color[rgb]{0,0,0}\put(-1124,-1261){\line( 0,-1){900}}
}%
{\color[rgb]{0,0,0}\put(-824,-1261){\line( 0,-1){900}}
}%
{\color[rgb]{0,0,0}\put(-524,-1261){\line( 0,-1){900}}
}%
{\color[rgb]{0,0,0}\put(-1124,-1261){\line( 1, 0){900}}
\put(-224,-1261){\line( 0,-1){900}}
}%
{\color[rgb]{0,0,0}\put(976,-1561){\line( 1, 0){900}}
}%
{\color[rgb]{0,0,0}\put(976,-1861){\line( 1, 0){900}}
}%
{\color[rgb]{0,0,0}\put(976,-2161){\line( 1, 0){900}}
}%
{\color[rgb]{0,0,0}\put(976,-1261){\line( 0,-1){900}}
}%
{\color[rgb]{0,0,0}\put(1276,-1261){\line( 0,-1){900}}
}%
{\color[rgb]{0,0,0}\put(1576,-1261){\line( 0,-1){900}}
}%
{\color[rgb]{0,0,0}\put(976,-1261){\line( 1, 0){900}}
\put(1876,-1261){\line( 0,-1){900}}
}%
{\color[rgb]{0,0,0}\put(3076,-1561){\line( 1, 0){900}}
}%
{\color[rgb]{0,0,0}\put(3076,-1861){\line( 1, 0){900}}
}%
{\color[rgb]{0,0,0}\put(3076,-2161){\line( 1, 0){900}}
}%
{\color[rgb]{0,0,0}\put(3076,-1261){\line( 0,-1){900}}
}%
{\color[rgb]{0,0,0}\put(3376,-1261){\line( 0,-1){900}}
}%
{\color[rgb]{0,0,0}\put(3676,-1261){\line( 0,-1){900}}
}%
{\color[rgb]{0,0,0}\put(3076,-1261){\line( 1, 0){900}}
\put(3976,-1261){\line( 0,-1){900}}
}%
\put(3826,-1486){\makebox(0,0)[b]{\smash{{\SetFigFont{12}{14.4}{\rmdefault}{\mddefault}{\updefault}{\color[rgb]{0,0,0}$2$}%
}}}}
\put(1126,-2086){\makebox(0,0)[b]{\smash{{\SetFigFont{12}{14.4}{\rmdefault}{\mddefault}{\updefault}{\color[rgb]{0,0,0}$1$}%
}}}}
\put(1426,-1786){\makebox(0,0)[b]{\smash{{\SetFigFont{12}{14.4}{\rmdefault}{\mddefault}{\updefault}{\color[rgb]{0,0,0}$1$}%
}}}}
\put(-674,-2536){\makebox(0,0)[b]{\smash{{\SetFigFont{12}{14.4}{\rmdefault}{\mddefault}{\updefault}{\color[rgb]{0,0,0}$2$}%
}}}}
\put( 76,-1786){\makebox(0,0)[b]{\smash{{\SetFigFont{12}{14.4}{\rmdefault}{\mddefault}{\updefault}{\color[rgb]{0,0,0}$2$}%
}}}}
\put(-974,-1486){\makebox(0,0)[b]{\smash{{\SetFigFont{12}{14.4}{\rmdefault}{\mddefault}{\updefault}{\color[rgb]{0,0,0}$2$}%
}}}}
\put(-674,-1486){\makebox(0,0)[b]{\smash{{\SetFigFont{12}{14.4}{\rmdefault}{\mddefault}{\updefault}{\color[rgb]{0,0,0}$1$}%
}}}}
\put(-374,-1486){\makebox(0,0)[b]{\smash{{\SetFigFont{12}{14.4}{\rmdefault}{\mddefault}{\updefault}{\color[rgb]{0,0,0}$1$}%
}}}}
\put(-674,-1786){\makebox(0,0)[b]{\smash{{\SetFigFont{12}{14.4}{\rmdefault}{\mddefault}{\updefault}{\color[rgb]{0,0,0}$1$}%
}}}}
\put(-974,-1786){\makebox(0,0)[b]{\smash{{\SetFigFont{12}{14.4}{\rmdefault}{\mddefault}{\updefault}{\color[rgb]{0,0,0}$1$}%
}}}}
\put(-374,-2536){\makebox(0,0)[b]{\smash{{\SetFigFont{12}{14.4}{\rmdefault}{\mddefault}{\updefault}{\color[rgb]{0,0,0}$1$}%
}}}}
\put(-974,-2536){\makebox(0,0)[b]{\smash{{\SetFigFont{12}{14.4}{\rmdefault}{\mddefault}{\updefault}{\color[rgb]{0,0,0}$3$}%
}}}}
\put( 76,-1486){\makebox(0,0)[b]{\smash{{\SetFigFont{12}{14.4}{\rmdefault}{\mddefault}{\updefault}{\color[rgb]{0,0,0}$4$}%
}}}}
\put(1126,-1486){\makebox(0,0)[b]{\smash{{\SetFigFont{12}{14.4}{\rmdefault}{\mddefault}{\updefault}{\color[rgb]{0,0,0}$4$}%
}}}}
\put(1426,-1486){\makebox(0,0)[b]{\smash{{\SetFigFont{12}{14.4}{\rmdefault}{\mddefault}{\updefault}{\color[rgb]{0,0,0}$3$}%
}}}}
\put(1726,-1486){\makebox(0,0)[b]{\smash{{\SetFigFont{12}{14.4}{\rmdefault}{\mddefault}{\updefault}{\color[rgb]{0,0,0}$3$}%
}}}}
\put(1126,-1786){\makebox(0,0)[b]{\smash{{\SetFigFont{12}{14.4}{\rmdefault}{\mddefault}{\updefault}{\color[rgb]{0,0,0}$2$}%
}}}}
\put(2176,-1486){\makebox(0,0)[b]{\smash{{\SetFigFont{12}{14.4}{\rmdefault}{\mddefault}{\updefault}{\color[rgb]{0,0,0}$10$}%
}}}}
\put(2176,-1786){\makebox(0,0)[b]{\smash{{\SetFigFont{12}{14.4}{\rmdefault}{\mddefault}{\updefault}{\color[rgb]{0,0,0}$3$}%
}}}}
\put(1126,-2536){\makebox(0,0)[b]{\smash{{\SetFigFont{12}{14.4}{\rmdefault}{\mddefault}{\updefault}{\color[rgb]{0,0,0}$7$}%
}}}}
\put(1726,-2536){\makebox(0,0)[b]{\smash{{\SetFigFont{12}{14.4}{\rmdefault}{\mddefault}{\updefault}{\color[rgb]{0,0,0}$3$}%
}}}}
\put(1426,-2536){\makebox(0,0)[b]{\smash{{\SetFigFont{12}{14.4}{\rmdefault}{\mddefault}{\updefault}{\color[rgb]{0,0,0}$4$}%
}}}}
\put(2176,-2086){\makebox(0,0)[b]{\smash{{\SetFigFont{12}{14.4}{\rmdefault}{\mddefault}{\updefault}{\color[rgb]{0,0,0}$1$}%
}}}}
\put(3226,-1486){\makebox(0,0)[b]{\smash{{\SetFigFont{12}{14.4}{\rmdefault}{\mddefault}{\updefault}{\color[rgb]{0,0,0}$4$}%
}}}}
\put(3526,-1486){\makebox(0,0)[b]{\smash{{\SetFigFont{12}{14.4}{\rmdefault}{\mddefault}{\updefault}{\color[rgb]{0,0,0}$3$}%
}}}}
\put(3526,-1786){\makebox(0,0)[b]{\smash{{\SetFigFont{12}{14.4}{\rmdefault}{\mddefault}{\updefault}{\color[rgb]{0,0,0}$2$}%
}}}}
\put(3226,-1786){\makebox(0,0)[b]{\smash{{\SetFigFont{12}{14.4}{\rmdefault}{\mddefault}{\updefault}{\color[rgb]{0,0,0}$3$}%
}}}}
\put(3226,-2086){\makebox(0,0)[b]{\smash{{\SetFigFont{12}{14.4}{\rmdefault}{\mddefault}{\updefault}{\color[rgb]{0,0,0}$3$}%
}}}}
\put(3526,-2086){\makebox(0,0)[b]{\smash{{\SetFigFont{12}{14.4}{\rmdefault}{\mddefault}{\updefault}{\color[rgb]{0,0,0}$2$}%
}}}}
\put(3826,-1786){\makebox(0,0)[b]{\smash{{\SetFigFont{12}{14.4}{\rmdefault}{\mddefault}{\updefault}{\color[rgb]{0,0,0}$1$}%
}}}}
\put(3826,-2086){\makebox(0,0)[b]{\smash{{\SetFigFont{12}{14.4}{\rmdefault}{\mddefault}{\updefault}{\color[rgb]{0,0,0}$1$}%
}}}}
\put(4276,-1486){\makebox(0,0)[b]{\smash{{\SetFigFont{12}{14.4}{\rmdefault}{\mddefault}{\updefault}{\color[rgb]{0,0,0}$9$}%
}}}}
\put(3226,-2536){\makebox(0,0)[b]{\smash{{\SetFigFont{12}{14.4}{\rmdefault}{\mddefault}{\updefault}{\color[rgb]{0,0,0}$10$}%
}}}}
\put(3526,-2536){\makebox(0,0)[b]{\smash{{\SetFigFont{12}{14.4}{\rmdefault}{\mddefault}{\updefault}{\color[rgb]{0,0,0}$7$}%
}}}}
\put(3826,-2536){\makebox(0,0)[b]{\smash{{\SetFigFont{12}{14.4}{\rmdefault}{\mddefault}{\updefault}{\color[rgb]{0,0,0}$4$}%
}}}}
\put(4276,-1786){\makebox(0,0)[b]{\smash{{\SetFigFont{12}{14.4}{\rmdefault}{\mddefault}{\updefault}{\color[rgb]{0,0,0}$6$}%
}}}}
\put(4276,-2086){\makebox(0,0)[b]{\smash{{\SetFigFont{12}{14.4}{\rmdefault}{\mddefault}{\updefault}{\color[rgb]{0,0,0}$6$}%
}}}}
\end{picture}%

\end{center}
show that the representations $[2,1^4]$, $[4,3^2,2,1^2]$,
$[4,3^3,2^3,1^2]$ appear in the respective tensor products: 
$[4,2]\otimes[3,2,1]$, $[7,4,3]\otimes[10,3,1]$,
$[9,6,6]\otimes[10,7,4]$. 


\subsection{Application to a branching rule}\label{sec{sec:Sl_Sp}}

Here we apply Theorem~\ref{th:ssgrpe} to the subgroup $G=\SP(2n)$ of $\h
G=\GL(2n)$. This subgroup is spherical of minimal rank, so that 
 Theorem~\ref{th:ssgrpe} applies for any $\h\lambda$ and $\h w$.

We define $G$ as the subgroup of $\GL(2n)$ which preserves the
alternate form given by the matrix:
$$I=\begin{pmatrix}
J &0&\ldots&0\\
0&J&\ddots&\vdots\\
\vdots&\ddots&\ddots&0\\
0&\ldots&0&J
\end{pmatrix}
$$
where $J=
 \begin{pmatrix}
 0&1\\
 -1&0
\end{pmatrix}
$.
 Then we
choose for $\h T$ the group of invertible diagonal
matrices, and for any $i\in\{1,\ldots,2n\}$, we denote by $\h\varepsilon_i$ the
usual character of $\h T$. 
Set $\hlambda=\hlambda_1\h\varepsilon_1+\cdots+\hlambda_{2n}\h\varepsilon_{2n}$.
The Weyl group $\h W$ is the
 symmetric group $\mathfrak S_{2n}$ and 
$\h w^{-1}\hlambda=\hlambda_{\h w(1)}\h\varepsilon_1+\cdots+\hlambda_{\h w(2n)}\h\varepsilon_{2n}$,
for $\hw\in \mathfrak S_{2n}$.

Set $T=G\cap \h T$ and define, for any $i\in\{1,\ldots,n\}$,  the restriction
$\varepsilon_i=\rho(\h\varepsilon_{2i-1})$. Then
$(\varepsilon_1,\ldots,\varepsilon_n)$ is a basis of
characters of  $T$ and we have:
$$\rho(\h w^{-1}\hlambda)=(\hlambda_{\h w(1)}-\hlambda_{\h w(2)},\hlambda_{\h
  w(3)}-\hlambda_{\h w(4)},\ldots,\hlambda_{\h w(2n-1)}-\hlambda_{\h w(2n)}).$$ 

The Weyl group $W$ acts on the characters of $T$ by permuting 
coordinates and by multiplying some coordinates by $-1$. So
$\overline{\rho(\h w^{-1}\hlambda})$
is obtained by arranging in a weak decreasing order the absolute values
$|\hlambda_{\h w(2i-1)}-\hlambda_{\h w(2i)}|$, for $i\in\{1,\ldots,n\}$.
We summarize this in the following

\bigskip
\noindent
{\bf Rule.} {\it 1. Consider a permutation $(\hlambda_{\h
    w(1)},\ldots,\hlambda_{\h w(2n)})$ of the coordinates of a dominant
    weight $\hlambda$ of $\GL(2n)$.\\
2. Order the $n$ absolute values $|\hlambda_{\h w(2i-1)}-\hlambda_{\h w(2i)}|$ to obtain a dominant weight $\mu$
    of $G$.

Then the multiplicity of $V_G(\mu)$ in $\h V_{\h G}(\hlambda)$ is non-zero. }
\bigskip

We believe that this rule cannot be easily deduce from
the combinatorial rules as those explicated in \cite{Sunda}.

\section{Tensor product decomposition}


The aim of this section is to prove Theorem~\ref{thi:tensorprod} stated in 
the introduction. We also give, at the end, two examples.

\bigskip
\begin{remark}
  \begin{enumerate}
  \item Condition~\ref{ass:long} of Theorem~\ref{thi:tensorprod} 
implies that $\langle v\mu,\alpha^\vee\rangle\geq 0$ and $\langle w\nu,\alpha^\vee\rangle\geq 0$.
  \item 
  Theorem~\ref{thi:tensorprod} asserts that the half line generated by
$(\lambda,\,\mu,\,\nu)$ is saturated in the Littlewood-Richardson semigroup.

Indeed, assume that  $\lambda=v\mu+w\nu+k\alpha$ with a rational number $k$ satisfies
$(-w_0\lambda+\mu+\nu)_{|Z(G)}=0$.
We obtain that $-w_0\lambda+\mu+\nu=(\lambda-w_0\lambda)+(\mu-v\mu)+(\nu-w\nu)+k\alpha$.
But, $\lambda-w_0\lambda$, $\mu-v\mu$ and $\nu-w\nu$ belong to the root lattice. It follows
that $k\alpha$ has to belong to the root lattice and so $k$ is an integer.
\end{enumerate}
\end{remark}

\bigskip
The strategy of the proof of Theorem~\ref{thi:tensorprod} is 
similar to that of Theorem~\ref{th:ssgrpe2}. 
So, we first prove adaptations of Proposition~\ref{prop:asymp} and of Lemma~\ref{lem:H0X0}.

\subsection{Asymptotic version}

To prove Proposition~\ref{prop:asymp}, we used Lemma~\ref{lem:Kos} mainly due to B.~Kostant; here, in order to prove Lemma~\ref{lem:asymp2} below, we will need to use the following strongly result of semi-stability mainly due to D.~Luna. 

\begin{lemma}\label{lem:Cssnonvide}
Consider the variety $Y=(G/B)^3$.
Let $\lambda,\,\mu$ and $\nu$ be three dominant weights of $T$.
Let $\beta$ be a root of $(G,T)$.
Denote by $S$ the neutral component of the Kernel of $\beta$ in $T$.
Let $(u,v,w)\in W^3$ and $C$ be the irreducible component of $Y^S$ 
containing $(uB/B,\,vB/B,\,wB/B)$.   
We assume that $u\Phi^+\cap v\Phi^+\cap w\Phi^+$ contains $\beta$.

The following are equivalent:
\begin{enumerate}
\item\label{Cssnonvide} $C$ contains semistable points with respect to 
$\Li_\lambda\boxtimes\Li_\mu\boxtimes\Li_\nu$;
\item \label{GITconeC}$\left\{
  \begin{array}{llll}
    (u\lambda+v\mu+w\nu)_{|S}=0\\
\langle u\lambda,\beta^\vee\rangle+\langle v\mu,\beta^\vee\rangle-\langle w\nu,\beta^\vee\rangle\geq 0,\\
\langle u\lambda,\beta^\vee\rangle-\langle v\mu,\beta^\vee\rangle+\langle w\nu,\beta^\vee\rangle\geq 0,\\
-\langle u\lambda,\beta^\vee\rangle+\langle v\mu,\beta^\vee\rangle+\langle w\nu,\beta^\vee\rangle\geq 0.\\
  \end{array}
\right.$
\end{enumerate}
\end{lemma}

\begin{proof}
Let $L$ be the centralizer of $S$ in $G$; it is a Levi subgroup of $G$ of semisimple rank one.
The variety $C$ is isomorphic to the product of three copies of the complete flag manifold of $L$, 
i.e. $(\PP^1)^3$.
Moreover, $(\Li_\lambda\boxtimes\Li_\mu\boxtimes\Li_\nu)_{|C}$ is isomorphic as an abstract line bundle to 
$\O(\langle u\lambda,\beta^\vee\rangle)\boxtimes\O(\langle v\mu,\beta^\vee\rangle)\boxtimes
\O(\langle w\nu,\beta^\vee\rangle)$. 
Note that $\langle u\lambda,\beta^\vee\rangle$, $\langle v\mu,\beta^\vee\rangle$ and $\langle w\nu,\beta^\vee\rangle$ are non-negative integers, because $\beta\in u\Phi^+\cap v\Phi^+\cap w\Phi^+$.

It is not difficult to check that $(\PP^1)^3$ has semistable points for the action of ${\rm SL}_2$ or ${\rm PSL}_2$ with respect to $\O(a)\boxtimes\O(b)\boxtimes\O(c)$ (where $a$, $b$ and $c$ are non-negative integers) if and only if we have $$\left\{
  \begin{array}{lll}
    
a+b-c\geq 0,\\
a-b+c\geq 0,\\
-a+b+c\geq 0.\\
  \end{array}
\right.$$

Now, the first equation of \ref{GITconeC} means that $S$ acts trivially on 
 $(\Li_\lambda\boxtimes\Li_\mu\boxtimes\Li_\nu)_{|C}$; and so, induces a $L/S$-linearized line bundle on $C$.
The three inequalities of \ref{GITconeC} are equivalent to the fact that $C$ contains
semistable points for the action of $L/S$ (which is isomorphic to ${\rm SL}_2$ or ${\rm PSL}_2$) 
with respect to $(\Li_\lambda\boxtimes\Li_\mu\boxtimes\Li_\nu)_{|C}$.
Now, it is clear that condition~\ref{Cssnonvide} implies condition~\ref{GITconeC}.

The converse implication is a direct application of \cite[Corollary~2 and Remark~1]{Luna:adh}
(see also \cite[Proposition~8]{GITEigen} for a formulation that can be directly
applied here).
\end{proof}

We use notation of Section~\ref{sec:ssgrpe} with $\hG=G\times G$. In particular,
$X^\circ_{v,w}$ is the $G$-orbit of $(vB/B,wB/B)$ in $X=(G/B)^2$.

We now prove the adaptation of Lemma~\ref{lem:asym}.

\begin{lemma}\label{lem:asymp2}
With assumptions of Theorem~\ref{thi:tensorprod},
there exist $n>0$ and a section $\tau$ of $(\Li_\mu\boxtimes\Li_\nu)^{\otimes n}$ 
of weight $-\lambda$ for $B^-$
whose restriction to $X_{v,s_\alpha w}$ is non-zero.  
\end{lemma}

\begin{proof}
  We apply Lemma~\ref{lem:Cssnonvide} with the dominant weights $-w_0\lambda$, $\mu$ and $\nu$, 
the root $\alpha$ and $(s_\alpha w_0,\,v,\,w)\in W^3$. 
Then, the first equation of condition~\ref{GITconeC} of Lemma~\ref{lem:Cssnonvide} is clearly satisfied and the three inequalities of condition~\ref{GITconeC} are respectively equivalent to $\left\{
  \begin{array}{lll}
k\leq\langle v\mu,\alpha^\vee\rangle,\\
k\leq\langle w\nu,\alpha^\vee\rangle,\\
k\geq 0.\\
  \end{array}
\right.$

We now remark, because of condition~\ref{ass:long}, that $\{w_0B/B\}\times X_{v,s_\alpha w}$ intersects 
$C\cap (\{w_0B/B\}\times X)$ along an open subset, and  we conclude the proof of the lemma, 
using the same arguments as in Lemma~\ref{lem:asym}.
\end{proof}
\subsection{Proof of Theorem~\ref{thi:tensorprod}}

In this section, we suppose that all  assumptions of  Theorem~\ref{thi:tensorprod}
are fulfilled.
We also set $\bar w=s_\alpha w$ and we denote by $G_{u,\bar w}$ the isotropy subgroup of $(vB/B,\bar wB/B)$ in $G$, i.e. $G_{v,\bar w}=vBv^{-1}\cap \bar wB\bar w^{-1}$.

We now prove the equivalent of Lemma~\ref{lem:H0X0}.

\begin{lemma}\label{lem:H0X02}
 The space $\CC[G]^{(B^-)_{-\lambda}\times (G_{v,\bar w})_{v\mu+\bar w\nu}}$
has dimension one.  
\end{lemma}

\begin{proof}
  We first prove that
$$
\CC[G]^{(B^-)_{-\lambda}\times (T)_{v\mu+\bar w\nu}}
$$
has dimension one.
Let us recall a classical property of some characters of the representation $V_G(\lambda)$:
the weights $\lambda-l\alpha$ with $l\in\{0,\cdots,\langle\lambda,\alpha^\vee\rangle\}$ have exactly multiplicity 
 one for $T$.
Frobenius' theorem implies that $\CC[G]^{(B^-)_{-\lambda}\times (T)_{v\mu+\bar w\nu}}$ is isomorphic
to $V_G(\lambda)^{(T)_{v\mu+\bar w\nu}}$. 
Assumption~\ref{ass:rellmn} of Theorem~\ref{thi:tensorprod} implies that $v\mu+\bar w\nu=
\lambda-(\langle w\nu,\alpha^\vee\rangle-k)\alpha$.
Assumption~\ref{ass:ineq} of the same theorem implies that 
$0\leq \langle w\nu,\alpha^\vee\rangle-k\leq \langle\lambda,\alpha^\vee\rangle$.
We obtain the dimension of $\CC[G]^{(B^-)_{-\lambda}\times (T)_{v\mu+\bar w\nu}}$
from the above-mentioned classical property.\\

Let $s$ be a non-zero element of 
$\CC[G]^{(B^-)_{-\lambda}\times (T)_{v\mu+\bar w\nu}}$.
Since $T$ is contained in $G_{v,\bar w}$, it is sufficient to prove that for any $h\in G_{v,\bar w}$, 
we have:
$$
hs=(v\mu+\bar w\nu)(h)s.
$$
By Lemma~\ref{lem:asymp2}, there exist $n$ and 
a non-zero 
$s_n\in \CC[G]^{(B^-)_{-n\lambda}\times (G_{v,\bar w})_{nv\mu+n\bar w\nu}}$.
Consider the algebra $\mathcal A=\oplus_{n\geq 0}\CC[G]^{(B^-)_{-n\lambda}}$.
Now, in $\mathcal A$, $s^{\otimes n}$ is a non-zero element of 
$\CC[G]^{(B^-)_{-n\lambda}\times (G_{v,\bar w})_{nv\mu+n\bar w\nu}}$.
By  the first part of the proof, $s^{\otimes n}$ and $s_n$ have to be proportional. 
It follows that for any $h$ in $G_{v,\bar w}$
$$
(hs)^{\otimes n}=h.s^{\otimes n}=(nv\mu+n\bar w\nu)(h)s^{\otimes n}=((v\mu+\bar w\nu)(h)s)^{\otimes n}.
$$
Since $\mathcal A$ is the algebra of regular sections of powers of an ample
 line bundle over a $P\backslash G$, it is integrally closed.
But, for any $h\in G_{v,\bar w}$, $(\frac{hs}{s})^{\otimes n}$ and $(\frac{s}{hs})^{\otimes n}$
belong to $\mathcal A$. So, $hs$ and $s$ are proportional.
There exists a regular map $\theta\,:\,H\longto \CC^*$ such that
$$
hs=\theta(h)s
$$ 
for any $h\in G_{v,\bar w}$.
We easily check that $\theta$ must be a character of $G_{v,\bar w}$. 
But, the restriction of $\theta$ to $T$ equals $v\mu+\bar w\nu$ so that $\theta=v\mu+\bar w\nu$. The lemma follows.
\end{proof}

We now prove Theorem~\ref{thi:tensorprod}.

\begin{proof}
It remains to prove that $V_G(\lambda)^*$ is a submodule of $V_G(\mu)^*\otimes V_G(\nu)^*$.
We interpret the latter module as the space of sections of $\Li_\mu\boxtimes\Li_\nu$ on $X$ and
consider the following sequence of morphisms:
$$\xymatrix{
{\rm H}^0((G/B)^2,\Li_\mu\boxtimes\Li_\nu)\ar@{->>}[d]\\
{\rm H}^0(X_{v,\bar w},\Li_\mu\boxtimes\Li_\nu)\ar@{^{(}->}[d]\\
{\rm H}^0(X_{v,\bar w}^\circ,\Li_\mu\boxtimes\Li_\nu)\ar@{^{(}->>}[d]\\
\CC[G]^{\{1\}\times (G_{v,\bar w})_{v\mu+\bar w\nu}}
}
$$

The surjectivity of the first map is a particular case (known before) of Theorem~\ref{th:Brion}.
The injectivity of the second map is obvious. And the next isomorphism is obtained  by applying Frobenius' theorem.\\

 Now, by Lemma~\ref{lem:H0X02}, there exists a non-zero section $\sigma$ of $\Li_\mu\boxtimes\Li_\nu$ on $X_{v,\bar w}^\circ$ of weight $-\lambda$ for $B^-$. Then, for some $n>0$, $\sigma^{\otimes n}$ extends to $X_{v,\bar w}$ by Lemmas~\ref{lem:asymp2} and \ref{lem:H0X02} together.  Since $X_{v,\bar w}$ is normal, it follows that $\sigma$ also extends to a regular section of $\Li_\mu\boxtimes\Li_\nu$ on $X_{v,\bar w}$. Thus, the theorem is proved.
\end{proof}

\subsection{Examples}

\paragraphe
In the following example, we will see that the hypothesis on $\alpha$ to be simple,  in Theorem~\ref{thi:tensorprod},  is necessary.
Consider $G=\Sp_4$. Denote by $\alpha_1$ and $\alpha_2$ respectively the short and the long simple roots, and $\omega_1$ and $\omega_2$ the associated fundamental weights. Let $\mu=\nu=\omega_2$ (and $v=w=Id$). Then we can compute that $$V_G(\mu)\otimes V_G(\nu)=V_G(0)\oplus V_G(2\omega_1)\oplus V_G(2\omega_2).$$ Define $\lambda:=v\mu+w\nu-(\alpha_1+\alpha_2)=\omega_2.$ Note that $\lambda$ satisfies the conditions~\ref{ass:rellmn} and \ref{ass:ineq} of Theorem~\ref{thi:tensorprod} with $\alpha=\alpha_1+\alpha_2$, because $\langle \omega_2,(\alpha_1+\alpha_2)^\vee\rangle=2$. 
We cannot apply Theorem~\ref{thi:tensorprod} just because $\alpha_1+\alpha_2$ is not a simple root. And in fact, $V_G(\lambda)$ is not a submodule of $V_G(\mu)\otimes V_G(\nu)$.\\

\paragraphe\label{sec:expletensor}
In that section, we look at the positions of the dominant weights $\lambda$ obtained 
in Theorem~\ref{thi:tensorprod} for fixed $\mu$, $\nu$, $\alpha$ and varying $k$.
We prove that, by this way, we obtain an ``integral segment'' with at least one extremity
corresponding to an original PRV component.

\begin{prop}\label{prop:segment}
Let $\lambda$ be  a dominant weight as in Theorem~\ref{thi:tensorprod}. 
Suppose, for convenience, that $\langle v\mu,\alpha^\vee\rangle\leq\langle w\nu,\alpha^\vee\rangle$.
Set $k_{max}=\langle v\mu,\alpha^\vee\rangle$ and  $\lambda_k= v\mu+w\nu-k\alpha$. 
Let $k_0$ be such that $\lambda=\lambda_{k_0}$.

Then, for any $k_0\leq k\leq k_{max}$, $\lambda_k$ is a dominant weight.
Moreover, $V_G(\lambda_{k_{max}})=
V_G(s_\alpha v\mu+w\nu)$ is an original PRV component of $V_G(\mu)\otimes V_G(\nu)$.  

\end{prop}

\begin{proof}
Denote by $S$ the set of simple roots of $(G,B)$ and by $\omega_\gamma$ the fundamental weight 
corresponding to the simple root $\gamma$. 
Then, for all $0\leq k\leq k_{max}$, we can write 
$\lambda_k=\sum_{\gamma\in S}a_{\gamma,k}\omega_\gamma$, with the $(a_{\gamma,k})$'s in $\ZZ$. 
Note that 
$$a_{\alpha,k_{max}}=\langle\lambda_{k_{max}},\alpha^\vee\rangle
=-\langle v\mu,\alpha^\vee\rangle+\langle w\nu,\alpha^\vee\rangle\geq 0.$$  
Remark also that 
$$\alpha=\sum_{\gamma\in S}\langle\alpha,\gamma^\vee\rangle\omega_\gamma=
\sum_{\gamma\in S}b_\gamma\omega_\gamma,\mbox{ with }b_\alpha=
2\mbox{ and }b_\gamma\leq 0,\,\forall\gamma\neq\alpha.$$ 
Then, $a_{\alpha,k}$ decreases when $k$ increases, and for any $\gamma\neq\alpha$, 
$a_{\gamma,k}$ increases with $k$. 
Moreover, since $a_{\alpha,k_{max}}\geq 0$,  $a_{\alpha,k}$ is non-negative 
for all $0\leq k\leq k_{max}$. 
This implies that, as soon as $\lambda_k$ is dominant, it stays dominant when $k$ increases up to $k_{max}$. 
Now, the proposition follows from the fact that $\lambda=\lambda_{k_0}$ is dominant. 
\end{proof}

We now illustrate this proposition by the following example. Consider $G=\SL_3$ with simple roots $\alpha_1$ and $\alpha_2$. 
Let $\mu=7\omega_1+2\omega_2$ and $\nu=\omega_1+3\omega_2$. 
Then the following picture represents the set of dominant weights $\lambda$ such that $V_G(\lambda)$ is a submodule of $V_G(\mu)\otimes V_G(\nu)$.
In this example, $\mu+\nu$ is an element of the root lattice so that all weights of 
$V_G(\mu)\otimes V_G(\nu)$ are in the root lattice. 
Then, in order to make the picture nicer, we only draw the root lattice instead of 
the weight lattice.

\begin{figure}[htbp]
\begin{center}

\input{figure1.pstex_t}
\label{figure:exampleSL3}
\end{center}
\end{figure}

\bibliographystyle{amsalpha}
\bibliography{prv_MPRarxiv}


\end{document}